\documentclass[11pt]{article}
\usepackage{epsfig}
\usepackage{verbatim}
\usepackage{multirow}
\usepackage[pdfmenubar=false, letterpaper=false, pdfpagelabels=true]{hyperref}

\begin{document}

\title{On the so called Boy or Girl Paradox}
\author{G.~D'Agostini \\
Universit\`a ``La Sapienza'' and INFN, Roma, Italia \\
{\small (giulio.dagostini@roma1.infn.it,
 \url{http://www.roma1.infn.it/~dagos})}
}

\date{}

\maketitle

\begin{abstract}
A quite old problem has  been recently revitalized 
by Leonard Mlodinow's book {\it The Drunkard's Walk},
where it is presented in a way that has definitely confused
several people, that wonder why
the prevalence of the 
name of one daughter 
among the population
should change the probability
that the other child is a girl too. 
I try here to discuss the problem 
from scratch, showing that the rarity of the name plays
no role, unless the strange assumption of two identical
names in the same family is taken into account. 
But also the name itself does not
matter. What is really important is `identification', meant
in an acceptation  broader than usual, 
in the sense that a child 
is characterized by a set of attributes 
 that make
him/her uniquely identifiable (`that one') 
inside a family. 
The important point of how the information is acquired 
is also commented, suggesting an explanation of why several
people tend to consider the informations ``at least one boy''
and ``a well defined boy'' (elder/youngest or of a given name)
 equivalent.

\end{abstract}

\vspace{2.9cm}
\section{Introduction}
A classical series of problems in elementary probability
theory is about the gender combinations
($m$-$m$, $m$-$f$, $f$-$m$ and  $f$-$f$)
in a family of two children. Being this an academic exercise
(in the bad sense of the term), usually one does not
attempt to assess how much one believes that these combinations
happen in a real family. This means that the well known male over
female birth asymmetry is neglected, as are neglected
gender correlations within a family, like those
induced by genetic factors, or by the possibility of
monovular twins.

Once the conditions are properly defined, the
usual questions, besides the trivial one of male/female, are
\begin{enumerate}
\item[Q$_1$)] What is the probability of two boys?
\item[Q$_2$)] What is the probability of two boys, if the eldest child
          is a boy?
\item[Q$_3$)] What is the probability of two boys, if at least 
          one child is a boy?
\end{enumerate}

These questions can be promptly answered looking at
contingency table 1 that lists the space of the four 
equiprobable elementary cases.
\begin{table}
\begin{center}
\begin{tabular}{c|cc|c}
    Eldest & \multicolumn{3}{c}{Youngest}  \\
           & $m$   &  $f$      & $m\cup f$ \\
\hline
 $m$      & $1/4$ &   $1/4$   &  $1/2$ \\
 $f$      & $1/4$ &   $1/4$   &  $1/2$ \\
\hline 
$m\cup f$  & $1/2$ &    $1/2$ &  $1$  
\end{tabular}
\caption{{\sl Table of equiprobable cases of the 
four possible sequences of child's gender. The symbol
`$\cup$' stands for `OR'.}}
\end{center}
\end{table}
\begin{enumerate}
\item[A$_1$)] The probability of two boys is 1/4, or 25\%, 
             since it is just the probability of each
             elementary event, that all
             together have to sum up to unity, or 100\%.
\item[A$_2$)] If the eldest child is a boy, the space of 
             possibilities is squeezed to the first row of the
             table. We remain with two equiprobable cases,
             each of which gets probability 1/2. In formulae:
             \begin{eqnarray}
             P(Em\cap Y\!m\,|\,Em,I_0) &=& \frac{P(Em\cap Y\!m\,|\,I_0)}
                                               {P(Em\,|\,I_0)}
             \label{eq:A2} \\ 
                                      & =&   \frac{1/4}{1/2} = \frac{1}{2}\,.
             \end{eqnarray} 
             [The symbol `$\cap$' stands for a logical `AND'; `$|$' stands for  `given',
              or `conditioned by'; `$Em$' and `$Y\!m$' are short forms
              for ``the eldest is male'' and ``the youngest is male'';
              the condition $I_0$ is the 
              {\it background status of information}
              under which the probabilities are evaluated, that includes 
              the simplifying hypotheses stated above; when 
              there is a further condition, like $Em$ and  $I_0$
              in the l.h.s. of Eq.~(\ref{eq:A2}), 
              they are both indicated after the conditional symbol `$|$', 
              separated by a comma.]
\item[A$_3$)] Finally,
              the information that there is at least one boy in the family
              reduces the space of possibilities
              to three equiprobable cases, of which only one is 
              that of our interest, thus getting 1/3 (the symbol `$\cup$'
              in the following formulae indicated a logical `OR'). 
              Formally
             \begin{eqnarray}
             P(Em\cap Y\!m\,|\,Em\cup Y\!m,I_0) &=& \frac{P[(Em\cap Y\!m)\cap(Em\cup Y\!m)\,|\,I_0]}
                                               {P(Em\cup Y\!m\,|\,I_0)} \\ 
                &=&  \frac{P(Em\cap Y\!m\,|\,I_0}
                                               {P(Em\cup Y\!m\,|\,I_0)}  \\
                                      & =&   \frac{1/4}{3/4} = \frac{1}{3}
             \end{eqnarray}  
\end{enumerate}
Obviously, the problem can been turned into probability of girl-girl
by symmetry. 

The `complication' 
(mainly induced confusion) comes when the information about
the name of one child is provided (`Florida' in 
{\it The Drunkard's Walk}~\cite{Drunkard}):
\begin{description}
\item[$Q_4$] What is the probability of two boys, if one 
of the children is called Mark?\footnote{This question can be turned
into ``what is the probability that Mark has a brother or a sister?''
and any normal and sane person might wonder about the sense
of this {\it madness}, as said a friend of mine with zero math skill, 
when he saw a first draft of this paper on my desk, because 
-- he explained --
``it is absolutely equally likely that he has either a brother or a sister''.}
\end{description}

\section{How the child name changes the probabilities}
At this point the question is how table 1 is changed by the
information that one child is known by gender and name 
(the latter usually implying the former). 

The way the problem is often solved is to assume that the same
name can given to two different children in the same family. 
Frankly, I never heard of this possibility before. 
And, anyhow, if such a strange behavior occurs
in some very rare cases, it seems to me of an importance
much lower than all other questions that
have been neglected (male/female asymmetry, genetic biases, etc.).
What I found annoying is that this peculiar solution is not 
reported as a mathematical curiosity 
(see e.g.  the Drunkard's Walk or the 
 Wikipedia page site dedicated to
the so called `paradox'~\cite{Wikipedia} -- a puzzle one is
unable to solve is not necessarily a paradox), but 
as it would be `the solution'. 

Nevertheless, let us first see what happens when this possibility is
allowed. 
(By the way, one might think of families with children 
coming from previous marriages, in which case identical names
might occur, but this possibility is excluded in the 
often implicit assumptions of this kind of puzzles, 
often formulated as ``a lady has two children, \ldots''.)

\subsection{Allowing identical names for two children of a family}\label{ss:id_names}
Just to stay close to the formulation of 
the problem in 
the recent disputes
that have triggered this paper, let us focus on girl probabilities,
assuming we know that one of the child is a girl of a given
name. Splitting the female category into `female of that given name'
($fN$) and `female with any other name' ($f\overline N$),
there are now three possible cases for each child,
$\{m, fN, f\overline N\}$, no longer equiprobable. 

Calling $r$ the fraction of girls owning that name 
in the population, we get the following probabilities, 
under the new background condition  $I_1$:
$P(m\,|\,I_1)=1/2$, $P(fN\,|\,I_1)=r\times 1/2=r/2$ and
 $P(f\overline N\,|\,I_1)=(1-r)/2$. 
The nine possibilities and their probabilities,
calculated using the product rule 
(justified by elder/youngest name independence), 
are reported in table 2.
\begin{table}
\begin{center}
\begin{tabular}{cc||ccc|cc}
 \multicolumn{2}{c||}{Eldest} & \multicolumn{4}{c}{Youngest} &  \\
\hline\hline
             &   & $m$   &  \multicolumn{2}{c|}{$f$}      & $m\cup f$ &\\
             &   &     &  $fN$ &  $f\overline N$          &  & \\
\hline
 $m$   &           & $1/4$ &   $\mathbf{r/4}$   &  $(1-r)/4$                      & $1/2$ & \\

\multirow{2}{*}{$f$} & $fN$           & $\mathbf{r/4}$ &   $\mathbf{r^2/4} $   &  $\mathbf{r (1-r)/4} $  & $r/2$ &\multirow{2}{*}{$1/2$}\\
  & $f\overline N$ &  $(1-r)/4$ &  $\mathbf{ r (1-r)/4}$ & $(1-r)^2/4$ & $(1-r)/2$  & \\
\hline 
$m\cup f$ &      & $1/2$             &    $r/2$          &   $(1-r)/2$ & 1 &\\
& & & \multicolumn{2}{c|}{$1/2$} & &  
\end{tabular}
\caption{{\sl Table of probabilities of the possible cases
assuming that 
eldest and youngest children can have the same name (see text).}}
\end{center}
\end{table}
From the table we can calculate the probability of both females,
if we know one  girl by name:
\begin{eqnarray}
 P[(Ef\cap Y\!f)\,|\,(EfN\cup Y\!fN)\ ,I_1] &=& 
   \frac{P[(Ef\cap Y\!f)\cap(EfN\cup Y\!fN)\,|\,I_1]}
        {P[(EfN\cup Y\!fN)\,|\,I_1]}\,. 
\end{eqnarray}
The denominator is given by the five elements emphasized in boldface
in table 2, whose probability sum up to $r-r^2/4$.
The numerator is given by the three elements
that have $m$ neither in the rows nor in the columns,
whose probabilities are $r^2/4$, $r(1-r)/4$ and $r(1-r)/4$,
adding up to $(2r - r^2)/4$. We get then 
\begin{eqnarray}
 P[(Ef\cap Y\!f)\,|\,(EfN\cup Y\!fN)\ ,I_1] 
                                      & =&   \frac{(2 r-r^2)/4}{r-r^2/4} \\
                                      & =&   \frac{1}{2}\left[\frac{1-r/2}{1-r/4}\right]
\label{eq:p_r_2}\\
                                      & \approx& \frac{1}{2}- \frac{r}{8}
\hspace{1.3cm}\mbox{(for }r\ll 1\mbox{)} 
\end{eqnarray}
As we can see, 
the probability does depend on $r$, but it 
tends rapidly to 1/2 for small values of $r$, as
also shown in table 3 for some numerical values of this parameter.\footnote{
The value $r=0.02=1/50$ is that used in Ref.~\cite{Wikipedia}, for which
the probabilities of table 2 acquire the following values
\begin{center}
\begin{tabular}{lll|l}
0.2500 & 0.0050 & 0.2450 & 0.5000 \\
0.0050 & 0.0001 & 0.0049 & 0.0100 \\
0.2450 & 0.0049 & 0.2401 & 0.4900 \\
\hline
0.5000 & 0.1000 & 0.4900 & 1.0000
\end{tabular}
\end{center}
from which we get the following table of {\it expected values}
in 10000 families
\begin{center}
\begin{tabular}{rrr|r}
2500 & 50 & 2450 & 5000 \\
  50 &  1 &   49 &  100 \\
2450 & 49 & 2401 & 4900 \\
\hline
5000 & 100 & 4900 & 10000
\end{tabular}
\end{center}
[By the way, I would like to point out that
quoting expected values is a way to state \ldots what we expect
in a probabilistic sense -- and probability theory teaches how
to calculate standard expectation uncertainties 
(the $\sigma$'s) - 
and has little to do with `frequentistic approach', 
since the probabilities have not been evaluated by past frequencies
(statistical data).]
}
\begin{table}[h]
\begin{center}
\begin{tabular}{l|c}
 \multicolumn{1}{c|}{$r$}   &  $P(\mbox{two girls}\,|\,fN,I_1)$  \\
\hline
0.3   &  0.45946 \\
0.2   &  0.47368 \\
0.1   &  0.48718 \\
0.02  &  0.49749 \\
0.01  &  0.49875 \\
0.001 &  0.49988 \\
0.0001&  0.49999
\end{tabular}
\caption{{\sl Probability of two girls in  family, 
if we know by name a daughter, calculated 
as a function of the prevalence of that name within
the girls of that population. [Note, just for mathematical curiosity,
that if $r=1$ (all girls have the same name), Eq.~(\ref{eq:p_r_2}) 
gives a  probability
of 1/3, thus recovering $Q_3$. 
In fact, in this case telling the name adds no more information 
to ``at least one is female''.]}}
\end{center}
\end{table}
\newpage
\subsection{Unique names of children within a family}\label{eq:unique_name}
Let us now see what happens if we require that, as it 
normally happens, children names are unique.
The central element of table 2 goes to zero, but the sums
along rows and columns have to be preserved\footnote{If no further
information is provided, 
the probability
that any child is a female with the special name $N$ has to be 
$r/2$, no matter if the child in question is the eldest or the youngest.
Similarly, the probability of girl with a name different from $N$
has to be $(1-r)/2$.} 
[for example the probability of $EfN\cap Y\!f\overline{N}$
becomes $r(1-r)/4+r^2/4$, that is the same as 
$r/2-r/4$, i.e. $r/4$]. 
The result is shown in table 4 
(we label the central value of the table by `-' to
remark that this case is impossible by assumption). [Note
that the impossibility of identical children names 
constrains $r$ to be smaller than $1/2$, well above any 
reasonable value. Remember also that the probabilities
of table 4
reflect the several simplifying
assumptions of the problem.]
\begin{table}[h]
\begin{center}
\begin{tabular}{cc||ccc|cc}
 \multicolumn{2}{c||}{Eldest} & \multicolumn{5}{c}{Youngest}  \\
\hline\hline
             &   & $m$   &  \multicolumn{2}{c|}{$f$}      & $m\cup f$ &\\
             &   &     &  $fN$ &  $f\overline N$          &  & \\
\hline
 $m$   &           & $1/4$ &   $\mathbf{r/4}$   &  $(1-r)/4$                      & $1/2$ & \\
\multirow{2}{*}{$f$} & $fN$ & $\mathbf{r/4}$ &  -  & $\mathbf{r/4} $ 
& $r/2$ &\multirow{2}{*}{$1/2$}\\
  & $f\overline N$ &  $(1-r)/4$ &  $\mathbf{r/4}$ & $(1-2r)/4$ & $(1-r)/2$  & \\
\hline 
$m\cup f$ &      & $1/2$             &    $r/2$          &   $(1-r)/2$ & 1 &\\
& & & \multicolumn{2}{c|}{$1/2$} & &  
\end{tabular}
\caption{{\sl Same as table 2, but not allowing the identical names
of the children.
}}
\end{center}
\end{table}
\newpage
Contrary to table 2, the four cases that involve $fN$ are now 
{\it equiprobable} (each with probability $r/4$). It follows that the probability that the other
child is a boy or a girl is 50\%, {\it independently} of the 
rarity of the name. In formulae (note the new background condition $I_2$):
\begin{eqnarray}
 P[(Ef\cap Y\!f)\,|\,(EfN\cup Y\!fN)\ ,I_2] &=& 
                                           \frac{P[(Ef\cap Y\!f)\cap(EfN\cup Y\!fN)\,|\,I_2]}
                                                 {P(EfN\cup Y\!fN\,|\,I_2)}
\label{eq:PEfYf|fN} \\ 
                                      & =&   \frac{2\times r/4}{4\times r/4} \\
                                      & =&   \frac{1}{2}\,.\label{eq:PEfYf|fN=1/2}
\end{eqnarray}
\section{Does the name really matter?}
At this point it is easy to understand that we could replace 
$fN$ in the table by $fID$, where `$ID$' stands now for
`uniquely identified within the family', thus getting table 5.
\begin{table}[b]
\begin{center}
\begin{tabular}{cc||ccc|cc}
 \multicolumn{2}{c||}{Eldest} & \multicolumn{5}{c}{Youngest}   \\
\hline\hline
             &   & $m$   &  \multicolumn{2}{c|}{$f$}      & $m\cup f$ &\\
             &   &     &  $fID$ &  $f\overline{ID}$          &  & \\
\hline
 $m$   &           & $1/4$ &   $\mathbf{r/4}$   &  $(1-r)/4$                      & $1/2$ & \\
\multirow{2}{*}{$f$} & $fID$ & $\mathbf{r/4}$ &  -  & $\mathbf{r/4} $ 
& $r/2$ &\multirow{2}{*}{$1/2$}\\
  & $f\overline{ID}$ &  $(1-r)/4$ &  $\mathbf{r/4}$ & $(1-2r)/4$ & $(1-r)/2$  & \\
\hline 
$m\cup f$ &      & $1/2$             &    $r/2$          &   $(1-r)/2$ & 1 &\\
& & & \multicolumn{2}{c|}{$1/2$} & &  
\end{tabular}
\caption{{\sl Same as table 4, but based on `identification' of a girl.}}
\end{center}
\end{table}
Think, for example to the following statements
\begin{itemize}
\item ``the secretary of the department X of hospital Y in Rome is 
        daughter of my aunt B who has also another child'';
\item ``the parents of the actress starring in the last movie I have seen
        have two children'';
\item ``the mother of that lady has got two children'';
\item and so on\ldots
\end{itemize}
In all these cases the probability that 
the female in question has a sister is 50\%, 
as everybody that is not fooled by probability theory will
promptly tell us (see footnote 1). 
It is not just a question of knowing her
name, or knowing that she is the eldest or the youngest 
(that's the reason we recover the answer to $Q_2$!).
{\it What matters is that this person is somehow
uniquely `identified' in the family}, where `identified' is within
quote marks because it is not requested we know her passport number,
but just that we are able to point to her as {\it that one}.

If this is not the case, and two children could correspond to 
the same description, then table 2 holds, assuming no correlation
between the descriptions (if one is blond, there is high change 
that the other is blond too, and so on). Therefore we recover it
as table 6,
\begin{table}
\begin{center}
\begin{tabular}{cc||ccc|cc}
 \multicolumn{2}{c||}{Eldest} & \multicolumn{5}{c}{Youngest}   \\
\hline\hline
             &   & $m$   &  \multicolumn{2}{c|}{$f$}      & $m\cup f$ &\\
             &   &     &  $f{\cal ID}$ &  $f\overline{{\cal ID}}$          &  & \\
\hline
 $m$   &           & $1/4$ &   $\mathbf{r/4}$   &  $(1-r)/4$                      & $1/2$ & \\

\multirow{2}{*}{$f$} & $f{\cal ID}$           & $\mathbf{r/4}$ &   $\mathbf{r^2/4} $   &  $\mathbf{r (1-r)/4} $  & $r/2$ &\multirow{2}{*}{$1/2$}\\
  & $f\overline{{\cal ID}}$ &  $(1-r)/4$ &  $\mathbf{ r (1-r)/4}$ & $(1-r)^2/4$ & $(1-r)/2$  & \\
\hline 
$m\cup f$ &      & $1/2$             &    $r/2$          &   $(1-r)/2$ & 1 &\\
& & & \multicolumn{2}{c|}{$1/2$} & &  
\end{tabular}
\caption{{\sl Same as table 2, but with reference to {\it non unique 
identification} (note the symbol `${\cal ID}$' instead of `$ID$'
of table 5)
rather then name.}}
\end{center}
\end{table}
 but in terms of {\it non unique identification} (`${\cal ID}$') rather
than of name. Now it makes sense. 
In fact, since we are referring here to 
`identification' in a loose sense, it might really occur
that two daughters correspond to the same description
(`goes to college', or `play tennis', and so on).
Finally, the name can be considered a generic identification, 
in order to include
the possibility of identical names in a family (for example
in the cases of second marriages).
\section{Some Bayesian flavor}
Someone asks me about {\it the} Bayesian solution of the problem
(because I am supposed to be {\it a Bayesian}). 
But, besides the clarification that 
``I am not a Bayesian''~\cite{MaxEnt98},
such a kind of `alternative' solution of the problem 
does not exist. The solution is already that provided
by Eq.~(\ref{eq:PEfYf|fN}), because `Bayesians' 
just make use of probability theory to state
the relative beliefs of several hypotheses given some
well stated assumptions. In particular, the so called
Bayes' rule for this problem is essentially 
Eq.~(\ref{eq:PEfYf|fN}), that can be possibly  written 
in other convenient forms using the rules of probability.

\subsection{Reconditioning the probability of an hypothesis on the
light of a new status of information}
To make the point clearer, and calling $A=Ef\cap Y\!f$ 
(``both children are female'') and
$B=EfN\cup Y\!fN$ (``one child is a female of a given particular name'') 
to simplify the notation, we can rewrite
Eq.~(\ref{eq:PEfYf|fN}) as
\begin{eqnarray}
P(A\,|\,B,I_2)  &=& \frac{P(A\cap B\,|\,I_2)}
                         {P(B\,|\,I_2)} \label{eq:P_A_B_I2}\\
                &=& \frac{P(B\,|\,A,I_2)\,P(A\,|\,I_2)}
                         {P(B\,|\,I_2)}\,.
\end{eqnarray}
The latter expression
shows explicitly how the probability of $A$ is 
updated, by the extra condition $B$, via the factor 
$P(B\,|\,A,I_2)/P(B\,|\,I_2)$, i.e.
\begin{eqnarray}
P(A\,|\,B,I_2) &=& \frac{P(B\,|\,A,I_2)} 
                         {P(B\,|\,I_2)} \times P(A\,|\,I_2)\,.
\end{eqnarray}
The three ingredients we need to evaluate $P(A\,|\,B,I_2)$
can be easily read from table 4,
\begin{eqnarray}
P(A\,|\,I_2)   &=& \frac{1}{4} \\
P(B\,|\,I_2)   &=& r \\
P(B\,|\,A,I_2) &=& 2 r \,,
\end{eqnarray}
from which we get 
\begin{eqnarray}
P(A\,|\,B,I_2) &=& \frac{2 r}{r}\times\frac{1}{4} = 2\times\frac{1}{4} = \frac{1}{2}\,,
\end{eqnarray}
recovering the result of section \ref{eq:unique_name}
(note that it must be so because we are strictly using the probabilities
of table 4).

\subsection{Updating the odds}
We can do it in a different way, comparing the probability
of ``two girls'' ($A$) with that of ``only one girl'' 
[let us indicate the latter hypothesis as 
$C = (Ef \cap Y\!m) \cup (Em \cap Y\!f)$]. 
The probability of $C$ conditioned by $B$, i.e.
$P(C\,|\,B,I_2)$, could be obtained in analogy to 
Eq.~(\ref{eq:P_A_B_I2}), reading $P(C\cap B\,|\,I_2)$
from table 4. But it can be more instructive to get it
 {\it the Bayesian way}, using the formula that 
shows how relative probabilities are updated by the 
{\it Bayes factor} to take into account the new 
piece of information
(this second approach has also the advantage of 
getting rid of $r$ since the very beginning):
\begin{eqnarray}
\underbrace{\frac{P(A\,|\,B,I_2)}{P(C\,|\,B,I_2)}}_{\mbox{`updated odds'}} 
&=& 
\underbrace{\frac{P(B\,|\,A,I_2)}{P(B\,|\,C,I_2)}
           }_{\begin{array}{c} \mbox{`updating factor'} \\ \mbox{({\it Bayes factor})} \end{array}}
 \times
\underbrace{\frac{P(A\,|\,I_2)}{P(C\,|\,I_2)}}_{\mbox{`initial odds'}}\,.
\end{eqnarray}
The initial probability of two girls is one half that of a single 
girl, i.e.
\begin{eqnarray}
\frac{P(A\,|\,I_2)}{P(C\,|\,I_2)} &=& \frac{1}{2}\,,
\end{eqnarray}
while the probability that there is a girl with a precise name 
is proportional to the number of girls in the family
(remember that the condition `$I_2$' does not allow the same name), 
namely
\begin{eqnarray}
\frac{P(B\,|\,A,I_2)}{P(B\,|\,C,I_2)}  &=& 2\,.
\end{eqnarray}
It follows
\begin{eqnarray}
\frac{P(A\,|\,B,I_2)}{P(C\,|\,B,I_2)} & = & 
2 \times \frac{1}{2} = 1\,:
\end{eqnarray}
the girl of which we know the name has equal probability
to have a sister ($A$) or a brother ($C$), that is
\begin{eqnarray}
P(A\,|\,B,I_2) = P(C\,|\,B,I_2) & = & \frac{1}{2}\,.
\end{eqnarray}

\section{Conclusions}
The probability that, knowing the name of 
one child in a family of two, the other one child
is of the same gender has nothing to do with the rarity
of the name, unless the crazy possibility of identical 
names in a family is assumed (and if somebody insists
that this can happen, he/she is invited to
calculate more realistic probabilities that take into
account male/female asymmetry and genetic correlations;
also the possibility of identical names of children coming
from previous marriages are implicitly excluded in this kind of
puzzles, that usually talk of ``a lady having two children\ldots'').

Moreover, what matters is not the knowledge of the name, 
but rather something that allows us to point to him/her
as `that one'. For this reason $Q_2$ and $Q_4$ have the same solution.

I would like to end with some comments on the 
last of the three text book questions reminded in the introduction.
It seems to me that the reason there is
quite a broad tendency to confuse $Q_3$ with $Q_4$
(or similar questions involving the child identification, including $Q_2$),
is that {\it in normal life the information about boy/girl is
acquired simultaneously with other attributes that make 
the identification unique} (``my daughter Claudia''). 
People do not express themselves as in math
textbooks, 
stating that ``I have two children, and at least one of them 
is a boy'', or 
``my children are not both boys''. We usually gain 
this information in an indirect way. For this reason
several people have some initial difficulty to 
grasp that ``that lady has Claudia and another child''
is not the same as ``that lady has two children, at least one being a girl''. 

Moreover, even if a mother says ``if I had two boys'', we may
understand from the context (already knowing she has two children) 
that she has two girls, because we perceived that she 
emphasized  `boys'
instead of `two' (in the latter case we could 
think she has already a boy). Instead, 
if she said ``if my children were both boys'', 
we usually understand that she is expressing this way
because she has a boy and a girl.  
Therefore, besides stereotyped
recreational puzzles, 
the evaluation of probabilities, 
in the sense of how much we have to rationally believe
the several hypotheses, can be not trivial. We
need to take  properly
into account  all contextual information
``when the bare facts wont'do''~\cite{Pearl}.
Indeed, in probability evaluations not only the `facts' play a role, 
but also 
the words, their sound and the expression 
of the person who says them, and (too 
often ignored) the question to which they reply~\cite{Pearl}.

\vspace{0.5cm}
It is a pleasure to thank Dino Esposito, 
Enrico Franco, Paolo Agnoli, Serena Cenatiempo 
and Stefano Testa for the several interactions
on the issues discussed here and related ones.
The paper has benefitted from comments by Dino and Paolo.

\vspace{1.5cm}

{\small
\section*{Appendix --- On the direct calculation of the elements of table 4}
The elements of table 4 have been evaluated from the condition
that the `central' one vanishes and that the marginal probabilities
have to be preserved (this means that, for example, the probability
that the younger is female with the special name $N$ is $r/2$
if no other information is provided, because $r/2$ is the assumed
probability that an individual of that population carries that name). 
Nevertheless, one might be interested to calculate the eight non 
vanishing terms in a direct way. But this calculation might reserve
surprises, as we shall see. 

\subsection*{First row and first column of the table (at least one boy)}
Although the elements that contain at least one boy are the easiest 
ones to be evaluated, let us get them with some pedantic  
detail, for didascalic purposes and in preparation
of the less obvious cases.
 In particular, we shall
rewrite $Em\cap Y\!fN$ as  $Em \cap Y\!f \cap Y\!N$ to remember 
that we require the eldest child to be a boy, the youngest to be
a girl and the name of the girl to be the particular one in which 
we are interested ($Y\!N$). Applying the `chain rule' we get
\begin{eqnarray*}
P(Em \cap  Y\!m\,|\,I_2) & = & P(Em\,|\,I_2)\times P(Y\!m\,|\,Em,I_2) = 
                              \frac{1}{2}\times\frac{1}{2} = \frac{1}{4} 
\nonumber \\
P(Em \cap Y\!fN\,|\,I_2) & = & P(Em \cap Y\!f \cap Y\!N\,|\,I_2) \\
                        & = &  P(Em\,|\,I_2)\times P(Y\!f\,|\,Em,I_2) 
                               \times  P(Y\!N\,|\,Y\!f,Em,I_2) \\               
                       & = &  \frac{1}{2}\times\frac{1}{2}\times r
                        =    \frac{r}{4}  \\
P(Em\cap Y\!f\overline N\,|\,I_2) & = & P(Em \cap Y\!f \cap Y\!\overline N\,|\,I_2)  \\
                       &=&  P(Em\,|\,I_2)\times P(Y\!f\,|\,Em,I_2)  
                               \times  P(Y\!\overline N\,|\,Y\!f,Em,I_2) \\ 
                       & = &  \frac{1}{2}\times\frac{1}{2}\times (1-r) 
                        =    \frac{1-r}{4}\,,
\end{eqnarray*}
recovering the first row of table 4. (Note that $\overline N$ stands for 
`a feminine name different from $N$' and not `any name but $N$'!)

Similarly, the second and third elements of the first column are
\begin{eqnarray*}
P(EfN\cap Y\!m\,|\,I_2) & = & P(Ef\cap EN \cap Y\!m\,|\,I_2) \\
                     & = &  P(Ef\,|\,I_2)\times P(EN\,|\,Ef,I_2) 
                               \times  P(Y\!m\,|\,Ef,EN,I_2) \\     
                     & = &  \frac{1}{2}\times r \times \frac{1}{2} 
                       =   \frac{r}{4}\\
P(Ef\overline N \cap Y\!m\,|\,I_2) & = & P(Ef\cap E\!\overline N \cap Y\!m\,|\,I_2) \\
                     & = &  P(Ef\,|\,I_2)\times P(E\!\overline N\,|\,Ef,I_2) 
                               \times  P(Y\!m\,|\,Ef,E\overline N,I_2) \\     
                     & = &  \frac{1}{2}\times (1-r) \times \frac{1}{2} 
                       =    \frac{1-r}{4}\,.
\end{eqnarray*}
The first column of table 4 is also recovered.

\subsection*{Probabilities of both females}
The probability that two girls have the same name is zero, that is because
\begin{eqnarray*}
P(EfN \cap Y\!fN\,|\,I_2) & = & P(Ef\cap EN \cap Y\!f \cap Y\!N\,|\,I_2) \\
                     & = &  P(Ef\,|\,I_2)\times P(EN\,|\,Ef,I_2) 
                            \times  P(Y\!f\,|\,Ef,EN,I_2) \\ 
                     && \times  P(Y\!N\,|\,Y\!f,Ef,EN,I_2) \\     
                     & = &  \frac{1}{2}\times r \times \frac{1}{2} \times 0 = 0\,.
\end{eqnarray*}
The third element of the second row is given by
\begin{eqnarray*}
P(EfN \cap Y\!f\overline N\,|\,I_2) & = & P(Ef\cap EN \cap Y\!f \cap Y\!\overline N\,|\,I_2) \\
                     & = &  P(Ef\,|\,I_2)\times P(EN\,|\,Ef,I_2) 
                            \times  P(Y\!f\,|\,Ef,EN,I_2) \\ 
                     && \times  P(Y\!\overline N\,|\,Y\!f,Ef,EN,I_2) \\     
                     & = &  \frac{1}{2}\times r \times \frac{1}{2} \times 1 = \frac{r}{4}\,.
\end{eqnarray*}
(The fourth factor of the r.h.s. of the last equation is 1 because,
once we know the eldest child has the name $N$, the youngest cannot have 
that name.) Also the second row of table 4 is recovered. 

The missing elements of the third row might present some pitfalls.  
Let us start from $Ef\overline N \cap Y\!fN$, that is 
$Ef\cap E\overline N \cap Y\!f \cap Y\!N$. 
It can be calculated  as 
\begin{eqnarray*}
P(Ef\overline N\cap Y\!fN\,|\,I_2) &=& P(Ef\cap E\overline N \cap Y\!f \cap Y\!N\,|\,I_2) \\
                     &=&  P(Y\!f \cap Y\!N \cap Ef\cap E\overline N \,|\,I_2) \\
                     & = & P(Y\!f\,|\,I_2)\times P(Y\!N\,|\,Y\!f,I_2) 
                            \times  P(Ef\,|\,Y\!f,Y\!N,I_2) \\ 
                     && \times  P(E\overline N\,|\,Ef,Y\!f,YN,I_2) \\     
                     & = &  \frac{1}{2}\times r \times \frac{1}{2} \times 1 = \frac{r}{4}\,.
\end{eqnarray*}
thus obtaining the same value of the third element of the second row, 
that is $P(EfN \cap Y\!f\overline N\,|\,I_2)$, as we expect by symmetry and
as it was in table 4. 

\subsubsection*{A pitfall}
But one would like to calculate $P(Ef\overline N\cap Y\!fN\,|\,I_2)$
`the other way around`, i.e. applying the chain rule starting 
from $P(Ef\overline N\,|\,I_2)$.
If one tries to proceed this way, there is 
{\it high chance to arrive to the following result}
\begin{eqnarray*}
P(Ef\overline N\cap Y\!fN\,|\,I_2) &=& P(Ef\cap E\overline N \cap Y\!f \cap Y\!N\,|\,I_2) \\
                     & = & P(Ef\,|\,I_2)\times P(E\overline N\,|\,Ef,I_2) 
                            \times  P(Y\!f\,|\,Ef,E\overline N,I_2) \\ 
                     && \times  P(Y\!N\,|\,Ef,E\overline N,Y\!f,I_2) \\     
                     & = &  \frac{1}{2}\times (1-r) \times \frac{1}{2} \times r \\
                     & = &  \mathbf{\frac{r\,(1-r)}{4} =\frac{r}{4} -  \frac{r^2}{4}}\,,
\end{eqnarray*}
{\it that differs from the value $r/4$ got previously}. 

In a similar way, one could be tempted to evaluate 
the probability that there are two girls, 
none of them carrying the name $N$, as 
\begin{eqnarray*}
P(Ef\overline N\cap Y\!f\overline N\,|\,I_2) &=& P(Ef\cap E\overline N \cap Y\!f \cap Y\!\overline N\,|\,I_2) \\
                     & = & P(Ef\,|\,I_2)\times P(E\overline N\,|\,Ef,I_2) 
                            \times  P(Y\!f\,|\,Ef,E\overline N,I_2) \\ 
                     && \times  P(Y\!\overline N\,|\,Ef,E\overline N,Y\!f,I_2) \\     
                     & = &  \frac{1}{2}\times (1-r) \times \frac{1}{2} \times (1-r) \\
                     & = &  \mathbf{\frac{(1-r)^2}{4} =\frac{(1-2r)}{4} +  \frac{r^2}{4}}\,,
\end{eqnarray*}
thus obtaining table 7,
\begin{table}
\begin{center}
\begin{tabular}{cc||ccc|cc}
 \multicolumn{2}{c||}{Eldest} & \multicolumn{5}{c}{Youngest}  \\
\hline\hline
             &   & $m$   &  \multicolumn{2}{c|}{$f$}      & $m\cup f$ &\\
             &   &     &  $fN$ &  $f\overline N$          &  & \\
\hline
 $m$   &           & $1/4$ &   $\mathbf{r/4}$   &  $(1-r)/4$                      & $1/2$ & \\
\multirow{2}{*}{$f$} & $fN$ & $\mathbf{r/4}$ &  -  & $\mathbf{r/4} $ 
& $r/2$ &\multirow{2}{*}{$1/2$}\\
  & $f\overline N$ &  $(1-r)/4$ &  $\mathbf{r(1-r)/4}$ & $(1-r)^2/4$ & $(1-r)/2$  & \\
\hline 
$m\cup f$ &      & $1/2$             &    $r/2-r^2/4$  &   $(1-r)/2+r^2/4$ & 1 &\\
& & & \multicolumn{2}{c|}{$1/2$} & &  
\end{tabular}
\caption{{\sl Same as table 4, but obtained by a {\bf wrong} reasoning that 
implicitly assumes that the condition $Ef\overline N$ does 
not change the probabilities of $Y\!fN$ and of $Y\!f\overline N$.
}}
\end{center}
\end{table}
that differs from table 4. Namely, 
 in the case of two females, it is now less probable that 
      the youngest girl has the particular name $N$. 
The probabilities differ by $r^2/4$, thus being negligible
for small $r$. 

One might think it is right so, because it reflects the
order of naming the children (``since the name $N$ cannot be 
given twice, the the eldest girl has a kind of first choice''). 
But on the other hand, we are dealing
here with knowledge (or ignorance), and therefore $P(EfN\,|\,I_2)$
and $P(Y\!fN\,|\,I_2)$ must be absolutely equal. It is just a question
of symmetry in reasoning in conditions of uncertainty. It doesn't
matter if we start thinking from the eldest or from the youngest child.
Stated in different words, from a probabilistic point of view 
{\it `eldest' and
`youngest' are mere labels}. 
The probability `matrix' {\it must be} symmetric. 

Moreover, it is curious to realize that table 7 produces a probability 
of two females that depends on $r$, as we saw in section \ref{ss:id_names}.
We get in fact 
\begin{eqnarray*}
 P[(Ef\cap Y\!f)\,|\,(EfN\cup Y\!fN)\ ,I_2] &=& \frac{r/4 +  r(1-r)/4}
                                                 {3\times r/4 + r(1-r)/4}\\
                                      & =&\mathbf{   \frac{1}{2}\left[\frac{1-r/2}{1-r/4}\right]}\,,
\end{eqnarray*}
exactly the {\it same result of section \ref{ss:id_names}}
[see Eq.~(\ref{eq:p_r_2})]. 
Therefore those who maintain that the probability of two girls, 
provided we know that one child is girl known by name,
does depend 
on the rarity of the name either assume that identical names
are possible inside the same family (a bizarre assumption), 
or have been caught
by this pitfall (a mistake in reasoning).

\subsection*{Conditional probabilities of female names}
The weak points of the previous evaluations, that lead
to table 7, with all its consequences, are the conditional probabilities
$P(Y\!N\,|\,Ef,E\overline N,Y\!f,I_2)$ and 
$P(Y\!\overline N\,|\,Ef,E\overline N,Y\!f,I_2)$, for which we assumed
{\it intuitively}
the values $r$ and $(1-r)$, as if the information that the eldest child
is a girl with a name different from $N$ did not change
the probability of the name of the other girl. 
This intuition, 
{\it roughly but not exactly correct}, is due to the fact 
that we tend to consider $r$ small (any modern population  
has a large amount of possible feminine names) such that 
the assumption that a girl has any name but the particular
one ($N$) does not change sizably the probability of the name 
of the other girl. This is the reason why the correct results are 
recovered for $r\rightarrow 0$, 
that was the hidden initial assumption!

But, strictly speaking, the $Ef\overline N$ and  $Y\!fN$ are not 
independent (in probability, or `stochasticly'), as 
are not independent  $Ef\overline N$ and  $Y\!f\overline N$:
the information that the eldest girl has a name different from $N$
{\it has to} increase the probability that the youngest 
girl is called $N$ (and has to decrease the probability that 
also the youngest girl has a name different from $N$). 
Comparing tables 4 and 7 we see that the effect goes this direction 
and has a size that decreases rapidly with $r$, 
going as $r^2$.

From this reasoning we get the following qualitative results:
\begin{eqnarray*}
P(Y\!N\,|\,Ef,E\overline N,Y\!f,I_2) &>& P(Y\!N\,|\,Y\!f,I_2) \\
P(Y\!\overline N\,|\,Ef,E\overline N,Y\!f,I_2) &<& P(Y\!\overline N\,|\,Y\!f,I_2) \,.
\end{eqnarray*}
The only problem is that it is not easy to evaluate these probabilities. 
But, fortunately, they can be calculated from the general rules of probability,
remembering that, as discussed above, the probabilities of 
$Ef\overline N\cap Y\!fN$ and of $Ef\overline N\cap Y\!f\overline N$
can be obtained in different ways. More precisely, the probability
of $Ef\overline N\cap Y\!fN$ can be calculated either directly, from
the easy $P(Y\!\overline N\,|\,Y\!f,Ef,EN,I_2)$, as we have done
just above in the appendix, or indirectly, requiring 
$P(Y\!fN\,|\,I_2) = r/2$, as it was done building up table 4
in section \ref{eq:unique_name}.
Instead, the last element of the table,  
$P(Ef\overline N\cap Y\!f\overline N\,|\,I_2)$, can 
by only calculated indirectly, either requiring 
$P(Y\!f\overline N\,|\,I_2) = (1-r)/2$, or
from the normalization rule, i.e. from the
constraint that all elements of the table
have to sum up to one.

Therefore, the conditional probabilities of interest can be
finally evaluated from the joint probabilities as
\begin{eqnarray*}
P(Y\!N\,|\,Ef,E\overline N,Y\!f,I_2) &=& \frac{P(Y\!N\cap Ef \cap E\overline N \cap Y\!f\,|\,I_2)}
                                            {P(Ef \cap E\overline N \cap Y\!f\,|\,I_2)} \\
                                   &=& \frac{P(Y\!N\cap Ef \cap E\overline N \cap Y\!f\,|\,I_2)}
                                            {P(Ef \cap E\overline N\,|\,I_2)\cdot 
                                             P(Y\!f\,|\,Ef \cap E\overline N,I_2)} \\ 
                                   &=& \frac{r/4}{(1-r)/2\times 1/2} = \frac{r}{1-r}\\ 
                                   &\approx& r  \hspace{4.0cm}\mbox{(for }r\ll 1\mbox{)} \\
P(Y\!\overline N\,|\,Ef,E\overline N,Y\!f,I_2) &=&  \frac{P(Y\!\overline N\cap Ef \cap E\overline N \cap Y\!f\,|\,I_2)}
                                            {P(Ef \cap E\overline N \cap Y\!f\,|\,I_2)} \\
                                   &=& \frac{P(Y\!\overline N\cap Ef \cap E\overline N \cap Y\!f\,|\,I_2)}
                                            {P(Ef \cap E\overline N\,|\,I_2)\cdot 
                                             P(Y\!f\,|\,Ef \cap E\overline N,I_2)} \\ 
                                   &=& \frac{(1-2r)/4}{(1-r)/2\times 1/2} = \frac{1-2r}{1-r}\\
                                   &\approx& 1 - r \hspace{3.5cm}\mbox{(for }r\ll 1\mbox{)} \,.
\end{eqnarray*}
[Obviously, the latter probability could have been calculated easier
as 
$P(Y\!\overline N\,|\,Ef,E\overline N,Y\!f,I_2) = 1 - P(Y\!N\,|\,Ef,E\overline N,Y\!f,I_2) = 1 - r/(1-r)$,
getting the same result.] 

Note that for very small values of $r$
we recover $P(Y\!N\,|\,Y\!f,I_2)=r$ and $P(Y\!\overline N\,|\,Y\!f,I_2)=1-r$,
 respectively.
That is, in this limit $Ef\overline N$ and  $Y\!fN$, as well as 
$Ef\overline N$ and  $Y\!f\overline N$, are approximately independent,
in agreement with our initial intuition.

Finally, we remind that the probabilities of the eldest girl name, 
conditioned by 
$Y\!f\overline N$, 
can be obtained by symmetry, i.e. 
$P(E\!N\,|\,Yf,Y\overline N,E\!f,I_2) = P(Y\!N\,|\,Ef,E\overline N,Y\!f,I_2) = r/(1-r)$
and 
$P(E\!\overline N\,|\,Yf,Y\overline N,E\!f,I_2) 
= P(Y\!\overline N\,|\,Ef,E\overline N,Y\!f,I_2) = (1-2r)/(1-r)$, and that 
all these expressions depend on the simplifying assumptions embedded in this 
kind of recreational puzzle.
}
\end{document}